\documentstyle[12pt]{amsart}
\def\F{\hfill $\Box$}

\def\gz{\ifmmode{Z\hskip -4.8pt Z}
    \else{\hbox{$Z\hskip -4.8pt Z$}}\fi}
\newtheorem{thm}{Theorem}[section]

\newtheorem{lemma}[thm]{Lemma}
\newtheorem{rmk}[thm]{Remark}
\newtheorem{defn}[thm]{Definition}
\newtheorem{cor}[thm]{Corollary}

\newtheorem{DT}[thm]{Definition and Theorem}
\title[Certain Applications of the Burnside rings]{Certain applications of the Burnside rings and ghost rings in the representation theory of\\ finite groups (II)}
\date{}
\author{K. K. Nwabueze}
\address{K. K. Nwabueze\\
Mathematical Sciences Research Institute\\
1000 Centennial Drive\\
Berkeley, CA 94720-5070, U.S.A.}
\email{nwabueze@@msri.org}

\author{F. Van Oystaeyen}
\address{
F. Van Oystaeyen\\
Department of Maths and Computer Science\\ 
University of Ant\-wer\-pen (UIA)\\ 
2610 Wilrijk, BELGIUM}
\email{voyst@@uia.ua.ac.be}

\thanks{Dr. Nwabueze's Research at MSRI is supported in part by the
University of Antwerpen (UIA) and the United States National Science
Foundation grant DMS-9022140.} 
\begin{document}
\maketitle

\begin{abstract}
Using the Burnside ring theoretic methods a new setting and a complete description of the Artin exponent $A(G)$ 
of finite $p$-groups was obtained in \cite{Nw}. In this paper, we compute $A(G)$ for any finite 
group $G$ -- hence providing the global version of \cite{Nw}.
\end{abstract}

\section{Introduction}
This is a continuation of the paper \cite{Nw}, where we  obtained a new setting
for the Artin exponent of a finite group and described these exponents for 
finite $p$-groups. In this note we show that the restriction made in \cite{Nw}
 that $G$ is a $p$-group can be remove. In otherwords the 
Artin exponent for any finite group $G$ depends on the value of a $p$-subgroup 
of $G.$ We claim that the Artin exponent $A(G\ ,\ {\cal U})$ of a finite $G$ 
is equal to $1$ if and only if $G$ is cyclic whereas $A(G\ ,\ {\cal U})$ 
difers by $p$ from the order of a $p$-subgroup of $G$ or is equal to the order 
of a $p$-subgroup of $G$ if $G$ is noncylic and theorem 4.2 dictates which. 
If the featuring $p$-subgroups of $G$ is Dihedral, Quoternion or semi-dihedral,
then $A(G\ ,\ {\cal U})$ is $2$ or $4$ and theorem 4.4 dictates which.  

In section 3 we show that we can reduce our problem essentially to the case 
where $G$ is a finite $p$-groups. For the convenience of the reader, we recall 
in section 2 some results in Burnside ring theory and the useful 
results in \cite{Nw}.  
Finally in section 4, we described the Artin exponent for finite group $G.$
\medskip

{\bf Some notation:} $G$ in this paper denotes a finite group. Let 
${1_{G}}$ denote the  trivial subgroup of $G;$ $p,$ a prime number; $<T>$, 
subgroup of $G$ generated by $T;\ $ $<T>_{N},$ the normal subgroup of $G$ 
generated by $T;\ $ $g_{h},$ the $g$ conjugate of $h;\ $ $[g\ ,\ h],$ the commutator of $g$ and $h\ -\ $ in particular 
$ [S\ ,\ T] := < [s\ ,\ t]\ |\ s \in S,\ t \in T>. $ Let $G^{sub/\sim}$
denote a set of representatives of the conjugacy classes of subgroups of
$G$ and $Z(G)$ denote the center of $G.$ Let 
$G_{s} := \{ g \in G\ |\ gs = s\ \forall\  s \in S \},$ the stabilizer subgroup
of $s.$

\section{Preliminaries}
Let $G$ be a finite group.
The set of isomorphism classes of finite left $G$-sets form a commutative
semi ring $B^{+}(G)$ with respect to disjoint union and cartesian product.
Using the Grothendieck ring construction on $B^{+}(G)$ one obtains the
corresponding Grothendieck ring $B(G)$ associated with $B^{+}(G).$ This
resulting ring is called the Burnside ring $B(G)$ of $G.$ So the
Burnside ring is the Grothendieck ring of finite $G$-sets which is
generated as an algebra over $Z$ by the isomorphism classes of finite left
$G$-sets $S_{1}\ S_{2}\ .\ .\ .\ $ subject to the relations 
$ S_{1} - S_{2} = 0 \ \mbox{ if }  S_{1} \cong S_{2}; \ \ $
$ S_{1} + S_{2} - (S_{1} \sqcup Y) = 0 \ \ $ and
$ S_{1} \cdot S_{2} - (S_{1} \times S_{2}) = 0. $
In otherwords, the elements of the Burnside ring are the virtual $G$-sets,
that is the formal difference $ S_{1} - S_{2} $ of isomorphism classes of
finite $G$-sets $S_{1},\ S_{2}.$ Furthermore one has that 
$ S_{1} - S_{2}\ =\ S'_{1} - S'_{2}\ \Longleftrightarrow\ S_{3} \sqcup S_{1}
\sqcup S'_{2}\ \cong\ S_{3} \sqcup S_{2} \sqcup S'_{1} $ where $S'_{1},\ S'_{2}
$ and $S_{3}$ are $G$-sets. This of course is equivalent with $S_{1} \sqcup S'_{2}\ \cong\ S_{2}
\sqcup S'_{1}.$ For any subgroup $ U \leq G$ of $G$ the map
$\chi_{U} :\ S\ \longrightarrow\ \#\{ s \in S\ \ |\ \ U \leq G_{s} \},$
from $S$ onto the number of $U$ invariant elements in $S,$ induces a 
homomorphism from $ B(G)$ into $Z$ -- the integers. Using standard arguments concerning the Burnside ring
the following statements are easily verified (see \cite{Dr}).
\medskip

(a) All homomorphism from $B(G)$ into $Z$ is of the form $\chi_{U}$ defined above. 
\medskip

(b) For any subgroup $U,\ V \leq G$ one has that $\chi_{U}\ =\ \chi_{V}$
if and only if $ U \sim V$ ($U$ is $G$-conjugate to $V$).
\medskip

(c) The product map $ \chi := \mathop{\Pi}\limits_{i=1}^{k} \chi_{i} :
B(G) \longrightarrow \mathop{\Pi}\limits_{i=1}^{k} Z$ of $k$ different
homomorphisms is injective and maps $B(G)$ onto a subring of finite index
$\mathop{\Pi}\limits_{i=1}^{k}(N_{G}(U_{i}\ :\ U_{i})$ of $\mathop{\Pi}\limits_{i=1}^{k} Z$ -- this way 
identifying $\mathop{\Pi}\limits_{i=1}^{k} Z$ with
the integral closure $\tilde{B}(G)$ of $B(G)$ in its total quotient
ring $ Q \otimes_{z} B(G) \cong \mathop{\Pi}\limits_{i=1}^{k} Q$ 
(see \cite{Dr},\cite{Di}). The product $\mathop{\Pi}\limits_{i=1}^{k} Z$ is 
called the {\it Ghost ring} of $G.$
\medskip

It is well known (see \cite{Dr}) that there are canonical set of
congruence relation where by one can characterize the image $\chi(B(G))$
as a subgroup of $\mathop{\Pi}\limits_{i=1}^{k} Z.$ More exactly, if 
$ U \unlhd V \leq G $ then for a $G$-set $X$ the Burnside lemma applied
with respect to the $V/U$-set $X^{U},$ implies that $ \mathop\sum\limits_{i=i}^{k} \chi_{<{\overline v}>} (X) \equiv 0 (mod\ (V:U)), $ where $<{\overline v}>$
denotes the subgroup of $V \leq G$ generated by the coset $ {\overline v} \subseteq V.$
\medskip

\begin{lemma}
Identifying $B(G)$ with its image in ${\tilde B}(G)$ with respect to the map
$\mathop\Pi\limits_{U \in G^{sub/\sim}} \chi_{U}$ one has that 
$\{ n \in Z\ |\ n{\tilde B}(G) \subseteq B(G)\}\ =\ |G|Z.$
\end{lemma}

Because for $U,\ V \leq G,$ one has $\chi_{U} = \chi_{V}$ if and only if $U$
and $V$ are $G$-conjugate and for $x,\ y \in B(G),$ one has $\chi_{U}(x) =
\chi_{U}(y)$ for all $U \leq G$ if and only if $x = y$ it follows that we
can identify each $x \in B(G)$ with the associated map (also denoted by $x$) 
$U \longrightarrow \chi_{U}(x),$ from the set of all subgroups of $G$ into 
$Z.$ Therefore we can consider $B(G)$ in a canonical way as a subring of the 
ghost ring $\tilde{B}(G).$
For more detail on the Burnside rings see \cite{Dr}. We recall the following simple facts (see any standard book on group theory). 

\begin{defn} 
Let $G$ be a $2$-group with generators $g$ and $h.$ We say that:\\{\bf(A)} $G$ is quaternion {\bf (Q)} if $g,$ $h$ have the relations
$$ g^{2^{n}} = 1,\ \  h^{2} = g^{2^{n-1}},\ \ hgh^{-1} = g^{-1}.$$
{\bf (B)} $G$ is  dihedral {\bf (D)} if $g,$ $h$ have the relations  
$$ g^{2^{n}} = 1,\ \  h^{2} = 1,\ \ hgh^{-1} = g^{-1}.$$ 
{\bf (C)} $G$ is  semi-dihedral {\bf (SD)} if $g,$ $h$ have the relations  
$$ g^{2^{n}} = 1,\ \  h^{2} = 1,\ \ hgh^{-1} = g^{-1 + 2^{n-1}}.$$ 
\end{defn}
\medskip

\begin{lemma}  
Let $G$ be a noncyclic abelian $p$-group $(p \ne 2)$ and $U \unlhd G,$  
$|U| = p.$ 
Then there exists a normal subgroup $ V \unlhd G $ of $G$ of order $p^{2},$ 
containing $U$ and isomorphic to  $ C_{p} \times C_{p}. $
\end{lemma}

\begin{lemma}
 Let $G$ be a finite $p$-group and\ $U$\ a subgroup of $G$ 
of order less than\ $p^{\beta}.$\ For $g \in G$ we have that $<g\ ,\ U>$ is 
cyclic of order $p^{\beta}$ if and only if $U \subseteq <g>.$
\end{lemma}

\begin{lemma}
Let $G$ be a finite $p$-group and let $U$ and $H$ be subgroups of $G$ such 
that $U$ is a cyclic normal subgroup with $ U \subseteq H.$  Let ${\cal{C}}(H)$
 denote the set of all cyclic subgroups $V$ of $H$ with $ U \leq V $ and 
$(V : U) = p,$ and let $ {\cal{C^{\prime}}}(H) = {\cal{C^{\prime}}}_{U}(H)$  
denote the subset of $ {\cal{C}}(H),$ consisting of those $ V \in {\cal{C}}(H)$
 which are normal in $H.$ We have:\\ 
(a.)\ $|{\cal{C}}(H)|\ \equiv\ |{\cal{C'}}(H)|\ \mbox{(mod } p).$\\
(b.)\ $ H^{'} = H^{'}_{U} := \{ h \in H\ |\ [h\ ,\ H] \subseteq U \}$ is a 
normal subgroup of $H.$\\ (c.)\ $ {\cal{C^{'}}}(H) = {\cal{C}}(H^{'})$\\ 
(d.) If $G^{'} \subseteq H \unlhd G$ then 
$|{\cal{C}}(H)|\ =\ |{\cal{C}}(G^{'})|(mod\ p). $ 
\end{lemma}

\begin{lemma}
With the above notation, one has that if $|U| = p \ne 2$ and $[G\ ,\ G] \subseteq U,$ then the following are eqivalent:\\
(i) $|{\cal{C}}(G)| =  1 $ \\ (ii) $ | {\cal{C}}(G) | \not\equiv 
0(\mbox{mod }\ p) $\\ (iii) $ G $ is cyclic
\end{lemma}

\begin{lemma}
 If $G$ is abelian, then the following are equivalent:\\ 
(i) $|{\cal{C}}(G)| =  1 $ \\ (ii) $ | {\cal{C}}(G) | \not\equiv 
0(\mbox{mod }\ p) $\\ (iii) $ G $ is cyclic
\end{lemma}

\begin{lemma}
If $|U| = p \ne 2,$ then the following are eqivalent:\\
(i) $|{\cal{C'}}(G)| =  1 $ \\ (ii) $ | {\cal{C'}}(G) | \not\equiv 
0(\mbox{mod }\ p) $\\ (iii) $ G' $ is cyclic
\end{lemma}

\begin{lemma} 
With the definitions and assumptions of (3.1), if $ G $ is a $ 2$-group, 
one has\\
(a.)\ $ |{\cal C}(H)| \equiv 1(2)  \Rightarrow  Z(G) \mbox  { is cyclic. } $\\
(b.)\ If $[G\ ,\ G] \subseteq U$\ and\ $|{\cal C}(H)|$\ is odd, then $G$ is 
cyclic or nonabelian of order $8.$\\
(c.)\ If $ H \leq G $ with $ [H\ ,\ H] \subseteq U \subseteq H $ and 
$| \{ V \in {\cal C}(H)\ |\ V \subseteq H \}| \equiv 1 (\mbox{mod}\ 2),$ then 
$H$ is cyclic or nonabelian of order $8.$ \\
(d.)\ If $|{\cal C}(H)| \equiv 1 (\mbox{mod}\ 2),$ then either $  H'$ is 
cyclic or $G$ is nonabelian of order $8.$\\
(e.)\ If $ |{\cal C}(H)| \equiv 1(\mbox{mod}\ 2), $ then $G$ contains a cyclic 
normal subgroup of index $2.$\\
\end{lemma}
\medskip

In \cite{Nw}, we characterized the elements in $\tilde{B}(G)$ which are in 
$B(G)$  in the following way:
\begin{DT}
 For ${\cal{U}}$ a family of cyclic subgroups of $G.$ Define 
$e_{{\cal{U}}}  \in \tilde{B}(G)$ by $e_{{\cal{U}}}(U) = 1$ if $U$ is 
contained in  ${\cal U}$ and $0$  otherwise. Then one has that 
$|G| \cdot e_{{\cal{U}}}(U) \in B(G).$  Furthermore define 
$A(G\ ,\ {\cal{U}}) := min(n \in N \mid n \cdot e_{{\cal{U}}} \in B(G)).$ 
The integer $A(G\ ,\ {\cal U})$ is said to be the {\it Artin exponent} of $G.$
One also has that $ A(G\ ,\ {\cal{U}}) $ divides $|G|$ (see \cite{Nw}).
\end{DT}

In \cite{Nw} we also gave a Burnside ring theoretic proof of the 
following well known results of Lam (see \cite{La}).

\begin{thm}
Let $G$ be a finite $p$-group. One has that $A(G\ ,\ {\cal U}) = 1$ if and 
only if $G$ is cyclic.
\end{thm}

\begin{thm}
 Let $G$ be a noncyclic finite $p$-group $(p \ne 2)$ of order $p^{\alpha}.$ 
Then one has $A(G\ ,\ {\cal U}) = p^{\alpha -1}. $ 
\end{thm}

\begin{thm}
If $p = 2 $ and $G$ is a noncyclic $2$-group of order $2^{\alpha}$ then 
$A(G) = 2^{\alpha -1},$ excluding the cases where $G$ is the Quaternion or 
Dihedral in which cases one has that $A(G\ ,\ {\cal U}) = 2.$ 
\end{thm}

These results provides the complete computation of $A(G\ ,\ {\cal U})$ where 
$G$ is a $p$-group.
\medskip

\section{Reduction to $ p$-Groups}

Let $p$ be a prime and consider
$Z_{p},$ the localization of $ Z $ at $ p, $ that is\\ 
$Z_{p}\ :=\ \{ a/b\ \ |\ \ a \in Z,\ b \in Z - pZ \}\ \subseteq\ Q.$ Put 
$B_{p}(G)\ :=\ Z_{p}\ \otimes_{z}\ B(G),$ and 
$ \tilde{B}_{p}(G)\ :=\ Z_{p}\ \otimes_{z}\ \tilde{B}(G).$
One has that $ Z_{p} \cdot \chi ( B(G) )\ \cong\ B_{p}(G)\ :=\ Z_{p}\ \otimes_{z}\ B(G)$ coincides with the 
subgroup of ${\tilde B}_{p}(G)$ for which the relations
$\mathop\sum\  n \cdot e_{\cal U}\ \equiv\ 0(mod\ (V_{p}\ :\ U)),$
with $V_{p}$ representing a subgroup of $G$ between $U$ and $V$ for which 
$V_{p}/U$ is the Sylow $ p$-subgroup of $V/U$ (see \cite{Dr}).
Therefore $e_{\cal U}$ is characterized as a subgroup of $ \mathop\Pi\limits_{i=1}^{k} Z $ by those relations, which runs 
through all the prime divisors of $|G|.$ Hence we have;

\begin{lemma}
An element $n \cdot e_{{\cal{U}}} \in \tilde{B}(G)$ is contained in $B(G)$ 
if and only if for every $U \unlhd V \leq G$ with $(V : U)$ a prime power one 
has  
$ \mathop{\sum}\limits_{vU \in V/U}n \cdot e_{{\cal{U}}}(<v\ ,\ U>)\ \equiv\ \\ 0 (\mbox{mod}\ (V : U)),$  where 
$e_{{\cal{U}}}(<v\ ,\ U>) := \# \{ vU \in V/U\ |\ <v\ ,\ U>\  \in\ {\cal{U}} \}.$  
\end{lemma}

\begin{cor}
For $ U \unlhd V \leq G,$ assume that $ (V : U) = p^{\alpha} $ for some 
prime $ p.$ Consider the decomposition $ U = U_{p} \times U_{p^{'}} $ where 
$|U_{p}| = p^{\beta}$ is a power of $ p $ and $ |U_{p^{'}}| $ is prime to $p.$
Let $ V_{p} $ denote the Sylow subgroup of $ V. $ Then one has\\
(a) $ (V : U) = p^{\alpha} $ implies that $ V = V_{p} \propto U_{p^{'}} $ and
$ V = V_{p} \times U_{p^{'}}$ if and only if $ U \leq Z(V),$ the center of $V.$ (Here $ \propto, $ resp. $ \times $ denotes the semidirect product, resp.
direct product).\\
(b) $ (V : U) = (V_{p} : U_{p}).$\\
(c) If\  $ c(U\ ,\ V) := \#\{ vU \in V/U\ \ |\ \ <v\ ,\ U> \mbox{ is cyclic } \} $ and \\
$ c(U_{p}\ ,\ V_{p}) := \#\{ vU_{p} \in V_{p}/U_{p}\ \ |\ \ <v\ ,\ U_{p}> \mbox { is cyclic } \} $ then $ c(U\ ,\ V)\ =\ c(U_{p}\ ,\ V_{p}) $ if $ V = V_{p} 
\times U_{p^{'}}, $ that is $ U \leq Z(V). $
\end{cor}

\begin{rmk}
 Observe that for\ $ vU \in V/U $ the group $ <v\ ,\ U> $ is cyclic only if\\ 
$ v \in C_{V}(U)\ :=\ \{ w \in V\ |\ wu = uw \mbox { for all } u \in U \}.$ 
This implies that $(V : C_{V}(U))$ must always divide the number $"n"$ 
in order for $(V\ :\ U)$ to divide $n \cdot {\bf c}(U\ ,\ V).$  So we have the 
following strategy for computing the Artin exponent.
For each $p$ subgroup $V$ of $G,$ consider the set of all cyclic normal 
subgroups $U$ of $V$ with $U$ contained in the centre of $V$ and compute the 
minimum of those positive integers $n$ for which $(V : U)$ divides 
$n \cdot c(U\ ,\ V).$
Finally, it is easy to see that if $U \leq  Z(V),$ for a $p$-group $V$ such 
that $ U $ cyclic one has that, 
$ {\bf c}(U\ ,\ V) = {\bf c}(U/U^{p}\ ,\ V/U^{p}), $ where $ U^{p} $ denotes 
the unique subgroup of $U$ of index $p.$ So we can restrict ourselves to the 
case where $|U| = p.$ 
\end{rmk}
We now summarize with the following,
\begin{thm}
Let $G$ be a finite group. Then $A(G\ ,\ {\cal U})\ =\ A(G_{p}\ ,\ {\cal U}),$
where $G_{p}$ is a $p$-subgroup of $G.$
\end{thm}

\section{Global Discussion}
For this discussion it suffices to consider the following cases;
\medskip

\underline{\bf Case 1:} $G$ is cyclic.
\medskip

The first thing to note is that the proof of Theorem(2.8) in \cite{Nw} does not depend 
on the fact that $G$ is a $p$-group. For the sake of explicitness we
state the following global version of the theorem. That is we remove the condition that $G$ is a $p$-group.

\begin{thm}  
Let $G$ be a finite group. For $ n \in N $ and a subgroup $U \leq G,$ let 
$x_{n} \in \tilde{B}(G)$ be defined by 

$$
x_{n} = 
\begin{cases}
n,&  \text{ if $ U $  is cyclic }\\
0,&  \text{ otherwise }
\end{cases}
$$
Then $A(G\ ,\ {\cal U}) = 1$ if and only if 
$G$ is cyclic.
\end{thm}

{\bf Proof:} We check by congruences.  Assume $U \unlhd V \leq G$ such that
$U$ cyclic of order $p$ and $U$ is contained in the center of $V$, where 
$|V| = p^{\alpha},$ ($p$ a prime). Then $ x \in B(G)$ implies 
 $\mathop\sum\limits_{vU \in V/U}x(<v\ ,\ U>) = n \cdot \#\{vU \in V/U\ |\ <v\ ,\ U> \mbox{ is cyclic}\}.$ To compute $A(G\ ,\ {\cal U}),$ we let
$V = G,$ that is $|G| = p^{\alpha}.$ These assumption remains valid for the
balance of this paper. Now if $n = 1,$ then one has that\\
$\#\{gU \in G/U\ |\ <g\ ,\ U> \mbox { is cyclic }\}\ \equiv \ 0(mod\ (G\ :\ U)).$
So\\ 
$(G\ :\ U)\ \ |\ \ \#\{gU \in G/U\ |\ <g\ ,\ U> \mbox{ is cyclic } \}.$ This implies\\ 
$|G|\ |\ \#\{ g \in G\ \ |\ \ <g\ ,\ U> \mbox { is cyclic } \}.$
 So for all $g \in G,$ one has that  $<g\ ,\ U>$ is cyclic. 
This implies that $G$ is cyclic since if $G$ is noncyclic then 
from (2.3) there exists a normal subgroup $ W\ \unlhd\ G$ such that 
$W$ contains $ U $ and $ W\ \cong\ C_{p} \times C_{p}\ \cong\ U \times C_{p}\ =\ U \times <g>$ for some $g \in G.$
That is $ <g\ ,\ U>\ =\ W$ for some $g \in G.$  So $G$ is cyclic.
\medskip

The converse is obvious.
\F
\bigskip

\underline{{\bf Case 2:}} $G$ is noncylic. 
\medskip

\begin{thm}
With the assumptions of Theorem 4.1 one has that if $G$ is noncyclic then
$$ A(G\ ,\ {\cal U}) := 
\begin{cases}
p^{\alpha},&  \text{if $G'\ :=\ \{ g \in G\ |\ [g\ ,\ U] \subseteq U \}$ is 
cyclic}\\ 
p^{\alpha - 1},&  \text{ if otherwise}.
\end{cases} $$ 
\end{thm}

{\bf Proof:} 
Since $G$ is noncyclic, then one obviously has that 
$|G|$ does not divide $\#\{g \in G\ |\ <g\ ,\ U> \mbox{ is cyclic} \}.$
If $G$ is abelian recall from lemma 2.3 that there exist a normal subgroup 
$W$ of order $p^{2}$ of the form\\ 
$\{ W \unlhd G\ |\ U\ \leq\ W \cong C_{p} \times C_{p}\}$ and indeed we have 
(see lemma 2.5 and 2.7) that\\
${\cal C}(W)\ :=\# \{C\ \leq\ W\ |\ C \mbox{ is cyclic}\ C \supseteq U \}  
\equiv 0(mod\ p).$  This implies $A(G\ ,\ {\cal U}) = p^{\alpha -1}.$
\medskip

If $G$ be nonabelian and $G'$ is noncyclic then repeating the above argument
(by setting $G = G',$ one has that $A(G\ ,\ {\cal U}) = p^{\alpha - 1}.$
If $G'$ is cyclic then it is easy to derive (from 2.8) that
$\#\{g \in G\ |\ <g\ ,\ U> \mbox{ is cyclic }\}\ \not\equiv\ 0 (mod\ p),$ 
in which case $A(G\ ,\ {\cal U})\ =\ p^{\alpha}.$                                \F
\medskip

We now compute $A(G\ ,\ {\cal U})$ for the special groups $Q,\ D,\ $ and $SD.$

\begin{cor}
Keeping the assumptions of (4.1) one has that if $G$ is any of $Q,\ D,\ $ or 
$SD$ then

$$ A(G\ ,\ {\cal U}) := 
\begin{cases}
4,&  \text{if $G'$ is cyclic}\\ 
2,&  \text{ if otherwise}.
\end{cases}
$$
\end{cor}

{\bf Proof:} Using the notation in section 3 one has here that if 
${\cal C}(G) \equiv 0(mod\ 2)$ then by a repetition of the above arguments
one has that $A(G\ ,\ {\cal U}) = 2^{\alpha - 1}.$ If ${\cal C}(G) \not\equiv 0(mod\ 2)$ then, from lemma 2.9, $G$ is cyclic or a nonabelian group of order 
$8.$ But by our assumption $G$ is noncyclic and so $G$ must be a nonabelian
group of order $8.$  It is now easy to see (by simple computation using lemma 2.9 repeatedly) that $A(G) = 4$ when $G'$ is cyclic and $A(G) = 2$ otherwise. 
\F

\pagestyle{plain}
{}
\medskip

\end{document}